\documentclass[11pt]{amsart}
\usepackage{xcolor}
\usepackage{graphicx}
\usepackage{latexsym}
\usepackage{amsfonts,amsmath,amssymb}
\newtheorem{theorem}{Theorem}[section]

\newtheorem{lemma}[theorem]{Lemma}

\usepackage{color}
\def\eps{\varepsilon} 

\def\a{\alpha}
\def\be{\beta}

\def\part{\partial}

\def\b1{\bold 1}

\newcommand{\beq}{\begin{equation}}
\newcommand{\eeq}{\end{equation}}

\theoremstyle{remark}

\numberwithin{equation}{section}
\date{\today}

\begin{document}

\title[Condorcet winner]{On likelihood of a Condorcet winner for uniformly random and independent voter preferences}

\author{Boris Pittel}
\address{Department of Mathematics, The Ohio State University, $231$ West $18$-th Avenue, Columbus, Ohio $43210-1175$ (USA)}
\email{bgp@math.ohio-state.edu}


\begin{abstract}
We study a mathematical model of a voting contest with $m$ voters and $n$ candidates, each voter ranking the candidates in order of preference, without ties. A Condorcet winner is a candidate who gets more than $m/2$ votes in
pairwise contest with every other candidate. An ``impartial culture'' setting is the case when each voter chooses his/her candidate preference list uniformly at random from all $n!$ preference lists, and does it independently of all other voters. 
For impartial culture case, Robert May and Lisa Sauermann showed that when $m=2k-1$ is fixed ($k=2$ and $k>2$ respectively), and $n$ grows indefinitely, the probability of a Condorcet winner is small, of order $n^{-(k-1)/k}$. 
We show that if $m, n\to\infty$ grow indefinitely and $m\gg n^4$, then, for every {\it fixed} $\ell>0$, the probability is at most of order 
$n^{-\ell} +n^2/m^{1/2}\to 0$, as $n\to\infty$.

\end{abstract}

\maketitle

\section{Condorcet winner} Suppose we have $m$ voters and $n$ possible candidates to vote for. Each voter ranks candidates in the order from the most to the least preferred. A Condorcet winner is a candidate $j\in [n]$ such that for every $j'\in [n]\setminus \{j\}$ the number of voters who prefer $j$ to $j'$ is {\it strictly\/} above $m/2$. There can be at most one Condercet winner. Suppose that the voters' preferences are random and, more specifically, uniform and independent of each other. 
This is known as an ``impartial culture'' (IC) case. Let $Q(m,n)$ denote the probability that there is a Condorset winner. How does $Q(m,n)$ 
depend asymptotically on $m,n$? The literature on this problem IC is substantial. Its main focus is the case when $m\to\infty$ (large electorate) and $n$ (candidates pool size) is arbitrary, but {\it fixed}.  See Niemi, Weisberg \cite{Nie}, Gehrlein \cite{Geh}, Krishnamoorthy,  Raghavachari \cite{KrR}, Tsetlin, Regenwetter, Grofman \cite{Tse}, Caizergues, Durand, Noy, de Panafieu, Ravelomanana \cite{Cai}, and many other references therein. 

As a rule, analysis depended on convergence of an auxiliary multinomial distribution to a limiting multidimensional normal distribution,
of dimension $n!$. 
In particular, Krishnamoorthy and Raghavachari proved in \cite{KrR} that, for $n$ {\it fixed\/},
\begin{equation}\label{0.9}
\lim_{m\to\infty}Q_{m,n}=Q_n:=n\int\limits_{-\infty}^{\infty}\!\!G^{n-1}(y)\frac{e^{-y^2}}{\sqrt{\pi}}\,dy,
\end{equation}
and that $Q_n=O(n^{-1/2}\log n)$ if $n\to\infty$. Recently de Panafieu \cite{Eli} used a generating function approach to prove a stronger 
bound $Q_n=O(n^{-1}\log^{3/2} n)$. (Back in 1971, May \cite{May} claimed, without proofs, similar but different results.)
 
In contrast,  May \cite{May} and recently Sauermann \cite{Sau} considered a diametrically opposite case, when $m$, the number of voters, is fixed, but $n$, the number of candidates, tends to $\infty$. 
%Let $Q_{m,n}$ denote the probability that there is a Condorcet winner.
Via a remarkably original argument, it was proved in \cite{Sau} that, for $k$ fixed and $n$ large, 
\begin{equation}\label{1}
Q_{2k-1,n}=C_k\cdot n^{-(k-1)/k}+O_k\Bigl(\tfrac{(\log n)^{1/k}}{n}\Bigr).
\end{equation}
Here $0<C_k<\infty$ is given by
\begin{equation}\label{2}
C_k=\int_{\bold x\ge \bold 0} \exp(-\sigma_{k,2k-1}(\bold x))\,d\bold x,\quad \bold x=\{x_1,\dots, x_{2k-1}\},
\end{equation}
and $\sigma_{k,2k-1}(\bold x)=\sum_{I\subset [2k-1],\,|I|=k}\prod_{i\in I} x_i$ is the usual $k$-th elementary symmetric polynomial of the variables $x_1,\dots x_{2k-1}$. %It was proved in \cite{Lisa} that $C_k\le [(2k-1)!]^2$. (Just by looking at the integral, it is not obvious at all that it converges.)  
Thus, the probability of a Condercet winner tends to zero at an exact rate $n^{-(k-1)/k}$, which for $k>2$ considerably strengthens a result  in \cite{May}.  
It was also pointed out in \cite {Sau}, that for $k=2$ the above claim, with a weaker remainder term, was proved in  
\cite{May}.

In this paper we consider a case complementary to the already studied cases, when both $m$ and $n$ tend to infinity. In that case, a reduction to the multinomial distribution becomes very problematic, if possible at all, since the dimension of the multinomial distribution grows super-exponentially with $n$. 

For analysis, we borrow from \cite{Sau} a key idea of representing the probability that a generic candidate is Condorcet winner as an integral over the $m$-dimensional unit cube by generating the random preference lists via an $m\times n$ matrix with entries that are independent, $[0,1]$-uniform random variables. The integrand is $(n-1)$-th power of probability that a sum of $m$  independent binomials 
(with expectations dependent on the integration point) exceeds $m/2$. To handle this dependency we apply a two-sided inequality for the integrand,  a special case of Esseen's inequality for the distribution of a sum of {\it non-identically\/} distributed independent terms, with finite third-order central moments. We prove

\begin{theorem}\label{newthm1} For $Q_n$ defined in \eqref{0.9}: 
%\[
%I_n=n\int_{-\infty}^{\infty} G^{n-1}(y)e^{-y^2/2}\,dy, \quad G(y)=\frac{1}{\sqrt{2\pi}}\int_{-y}^{\infty} e^{-z^2/2}\,dz.\]
(a) If $m\gg n^4$, then $Q_{m,n}=Q_n+O(n^2m^{-1/2})$. (b) If $n\to\infty$, then, for every {\it fixed} $\ell> 0$, $Q_n
=O(n^{-\ell})$; in words, $Q_n$ tends to zero super-polynomially fast. So, 
$Q_{m,n}\to 0$, if $m\gg n^4$ and $n\to\infty$.
\end{theorem}
\noindent
The theorem sheds light on the Condorcet probability when $n$, the number of candidates, tends to infinity as well, but considerably slower than $m$, which is the case of a very large electorate and a moderately large pool of candidates. (Can the condition $m\gg n^4$ be relaxed?)
Our bound for $Q_n$ considerably strengthens the bounds obtained by Krishnamoorthy and Raghavachari, and by de Panafieu, that we cited above.

In conclusion, we should mention that Sauermann's approach was very recently used in \cite{Cai} for a more general case of an $\boldsymbol{\a}$-winner, i. e. a generic candidate  preferred to each other candidate $j$ by a fraction $\a_j$ of the voters, for the case of large $m$ and {\it fixed\/} $n$. The analysis was based on techniques from analytic combinatorics developed by Flajolet, Sedgewick \cite{Fla} and Pemantle,
Wilson \cite{Pem}. In pursuit of generality, the authors did not dwell on any special fractions $\boldsymbol \a$'s.
%The analysis was based on techniques from analytic combinatorics developed by Flajolet and Sedgewick \cite{Fla}, and for several variables  by Pemantle and Wilson \cite{Pem}.

\subsection{Proofs}  Following \cite{Sau}, we generate the random instance of of
impartial culture preferences via an $m\times n$ matrix $X=\{X_{a,b}\}$ with the entries $X_{a,b}$ independent, $[0,1]$-Uniforms. Each voter $a\in [m]$ ranks the $n$ candidates from most preferred to least preferred in the increasing order of the entries of the $a$-th row $\{X_{a,1},\dots, X_{a,n}\}$. Let us evaluate $Q_{m,n}$, the probability that candidate $1$ is a Condercet winner. The candidate $1$ is preferred by strictly more than $m/2$ voters to a candidate $j>1$ if and only if $\sum_{a\in [m]}\Bbb I(X_{a,1}< X_{a,j})>m/2$, where $\Bbb I(A)$ denote the indicator of an event $A$. The $(n-1)$ events in question are not independent, but they are independent, if {\it conditioned\/} on the event $A:=\{X_{1,1}=x_1,\dots, X_{m,1}=x_m\}$. Now, the $(n-1)$ sums 
\begin{equation}\label{n1}
N_j(\bold x):=\sum_{a\in [m]}\Bbb I(x_a< X_{a,j}),\quad 2\le j\le n,
\end{equation}
are {\it equidistributed independent sums of independent terms\/}. Therefore,  introducing $C_{m}=[0,1]^{m}$,
the probability of a Condorcet winner is given by
\begin{equation}\label{n1.1}
Q_{m,n}=n\int_{\bold x\in C_{m}}\Bbb P^{\,n-1}(N_2(\bold x)>m/2)\,d\bold x,
\end{equation}
as there can be at most one Condorcet winner. Probabilistically,
\[
Q_{m,n}=n\Bbb E\bigl[\Bbb P^{n-1}(N_2(\bold x)>m/2)\big|_{\bold x=\bold X}\bigr],
\]
where $\bold X=\{X_a\}_{a\in [m]}$, $X_a$ being independent $[0,1]$-Uniforms. 

\noindent{\bf Note.\/} Of course, we could have chosen any other continuous distribution instead of the $[0,1]$-Uniform. However the uniform distribution is perfect for the asymptotic analysis. Here is one immediate benefit of our choice: $\tilde X_{i,j}:=1-X_{i,j}$ are also $[0,1]$-Uniforms, and $\tilde X_{a,b}< \tilde X_{a',b'}$ if and only if $X_{a,b}>X_{a',b'}$. This implies that $Q_{m,n}$ is also the probability of a Condorcet ``loser'', i. e. a candidate who gets at most $m/2$ votes in a pairwise contest with every other candidate. 

For certainty, we consider the case $m=2k-1$. In this case, Condorcet winner is a candidate who gets at least $k$ votes in every pairwise contest.

To analyze an asymptotic behavior of the integral \eqref{n1.1}, we will use the following classic theorem by Esseen \cite{Ess}, cf. Feller \cite{Fel}. Let $Y_1,\dots, Y_{2k-1}$ be independent random variables such that $\Bbb E[Y_a]=0$, $\Bbb E[Y_a^2]<\infty$, $\Bbb E[|Y_a|^3]<\infty$. Denote $s_{2k-1}^2=\sum_a \Bbb E[Y_a^2]$, $r_{2k-1}=\sum_a \Bbb E[|Y_a|^3]$. Then 
\begin{align*}
\left|\Bbb P\left(\frac{\sum_a Y_a}{s_{2k-1}}\ge y\right)-G(y)\right|
\le \frac{6t_{2k-1}}{s_{2k-1}^3},\quad G(y):=\frac{1}{\sqrt{2\pi}}\int_y^{\infty} e^{-z^2/2}\,dz.\\
\end{align*}
(So, $\Phi(y):=1-G(y)$ is the probability distribution of the Gaussian random variable with zero mean and unit variance.) In the case of $N_j(\bold x)$, we have $Y_a=Y_a(\bold x)=\Bbb I(x_a< X_{a,j})-(1-x_a)$, with $\Bbb E[Y_a]=0$,
\begin{align*}
\Bbb E[Y_a^2]&=x_a(1-x_a),\\
\Bbb E[|Y_a|^3]&=(1-x_a)x_a^3+x_a(1-x_a)^3\le x_a(1-x_a).
\end{align*}
So, $\sigma^2_{2k-1}(\bold x)=\sum_a x_a(1-x_a)$, and 
\[
\frac{t_{2k-1}(\bold x)}{s_{2k-1}^3(\bold x)}\le \frac{\sigma^2_{2k-1}(\bold x)}{\sigma_{2k-1}^3(\bold x)}=\frac{1}
{\sigma_{2k-1}(\bold x)}.
\]
Therefore, denoting $y(\bold x):=\frac{k-\sum_a (1-x_a)}{\sigma_{2k-1}(\bold x)}$,
\begin{multline*}
\Bbb P(N_2(\bold x)\ge k)=\Bbb P\left(\frac{N_2(\bold x)-\sum_a (1-x_a)}{\sigma_{2k-1}(\bold x)}\ge
y(\bold x)\right)\\
= G(y(\bold x))+O(\sigma_{2k-1}^{-1}(\bold x)),
\end{multline*}
uniformly over $\bold x\in C_{2k-1}$. Consequently, for a fixed constant $\a$,
\begin{multline*}
\Bbb P^{\,n-1}(N_2(\bold x)\ge k)\!=G^{n-1}(y(\bold x))\\
+O\!\left(\sum_{\nu=1}^{n-1}\binom{n-1}{\nu}G^{n-\nu}(y(\bold x)) \a^{\nu}\sigma^{-\nu}_{2k-1}(\bold x)\!\right).
\end{multline*}
Neglecting $G^{n-\nu}y(\bold x)$, and integrating over $C_{2k-1}$, we get
\begin{multline}\label{n1.11}
n^{-1}Q_{2k-1,n}=\Bbb E[G^{n-1}(y(\bold X))]\\
+O\!\left(\sum_{\nu=1}^{n-1}\binom{n-1}{\nu}\a^{\nu}\,\Bbb E\bigl[\sigma^{-\nu}_{2k-1}(\bold X)
\bigr]\!\right).
\end{multline}
%Here, by H\"older inequality, for $p,q>1$ and $1/p+1q=1$,
%\[
%\Bbb E\bigl[G^{n-\nu}(y(\bold X)) \sigma^{-\nu}_{2k-1}(\bold X)\bigr]\le \Bbb E^{1/p}\bigl[G^{p(n-\nu)}(y(\bold X))\bigr]
%\cdot \Bbb E^{1/q}\bigl[\sigma^{-q\nu}_{2k-1}(\bold X)\bigr].
%\]
\begin{lemma}\label{nlem1} For $n\ll k^{1/2}$, the big-O term in \eqref{n1.11} is $O(nk^{-1/2})$.
%we have
%\begin{equation*}
%\Bbb E\bigl[\sigma_{2k-1}^{-\nu}(\bold X)\bigr] =O((c/k)^{\nu/2}),\quad c>e.
%\end{equation*}
\end{lemma}
\begin{proof}
Recall that $\sigma_{2k-1}(\bold x)=\left(\sum_a x_a(1-x_a)\right)^{1/2}$, thus it is zero at all $2^{2k-1}$ corners of $C_{2k-1}$. For $A\subseteq [2k-1]$, let $C_{2k-1}(A)\!=\!\{\bold x\in C_{2k-1}\!: \forall a\in A, x_a\le 1/2\}$. For $\bold x\in C_{2k-1}(A)$, we have
\[
\sum_a x_a(1-x_a)\ge \frac{1}{2}\sum_{a\in A} x_a+\frac{1}{2}\sum_{a\in A^c}(1-x_a).
\]
So,
\begin{multline*}
\int\limits_{\bold x\in C_{2k-1}}\!\!\!\!\!\!\sigma^{-\nu}_{2k-1}(\bold x)\,d\bold x\le 2^{\nu/2}\!\!\!\!\!\sum_{A\subseteq [2k-1]}\,\int\limits_{\bold x\in C_{2k-1}(A)}
\!\!\!\left(\sum_{a\in A}x_a +\sum_{a\in A^c}(1-x_a)\right)^{-\nu/2}\!\!\!\!\!\!\! d\bold x.\\
\end{multline*}
Switching to $\xi_a=2x_a$ for $a\in A$ and $\xi_a=2(1-x_a)$ for $a\in A^c$, and denoting $s(\boldsymbol\xi)=\sum_a \xi_a$, we obtain 
\begin{multline}\label{n1.2}
\int\limits_{\bold x\in C_{2k-1}}\!\!\!\!\sigma^{-\nu}_{2k-1}(\bold x)\,d\bold x\le 2^{\nu/2}\!\!\sum_{A\subseteq [2k-1]}
\!\!2^{-(2k-1)}\int\limits_{\boldsymbol\xi\in C_{2k-1}}\!\!\!\!\!s(\boldsymbol\xi)^{-\nu/2}\,d\boldsymbol\xi\\
=2^{\nu/2}\!\!\!\int\limits_{\boldsymbol\xi\in C_{2k-1}}\!\!\!\!s(\boldsymbol\xi)^{-\nu/2}\,d\boldsymbol\xi,
\end{multline}
as the number of subsets of $[2k-1]$ is $2^{2k-1}$. It remains to upper bound the last integral. 

The integral equals $\Bbb E[S_{2k-1}^{-\nu}]$, where $S_{2k-1}=\sum_{a\in 2k-1}X_a$. Consider first 
$C(1)=\{\boldsymbol\xi\in C_{2k-1}: s(\boldsymbol\xi)\le \be k\}$, $\be$ to be chosen later. Since the density of 
$S_{2k-1}$ is below $s^{2k-2}/(2k-2)!$, we see that
\begin{multline*}
\int\limits_{\boldsymbol\xi\in C(1)}\!\!\!\!s(\boldsymbol\xi)^{-\nu/2}\,d\boldsymbol\xi
\le\frac{1}{(2k-2)!}\int_0^{\be k} s^{2k-2-\nu/2}\,ds\\
=\frac{(\be k)^{2k-1-\nu/2}}{(2k-2)!\cdot(2k-1-\nu/2)}=O\left((\be k)^{-\nu/2}\left(\frac{e\be}{2}\right)^{2k}
\right),\\
\end{multline*}
since $\nu\ll k$ and $(2k-1)!=\Theta\bigl(k^{1/2}((2k-1)/e)^{2k-1}\bigr)$. As
\[
\int\limits_{\boldsymbol\xi\in C_{2k-1}\setminus C(1)}\!\!\!\!s(\boldsymbol\xi)^{-\nu/2}\,d\boldsymbol\xi
\le  (\be k)^{-\nu/2},
\]
we see that the integral over the whole $C_{2k-1}$ is also of order at most $(\be k)^{-\nu/2}$ if we choose $\be<2/e$.
Thus, the equation \eqref{n1.2} yields 
\begin{equation}\label{n1.31}
\Bbb E \bigl[\sigma_{2k-1}^{-\nu}(\bold X)\bigr]=\int\limits_{\bold x\in C_{2k-1}}\sigma^{-\nu}_{2k-1}(\bold x)\,d\bold x=O\Bigl(\left(\frac{2}{\be k}\right)^{\nu/2}\Bigr).
\end{equation}
Consequently,
\begin{multline*}
\sum_{\nu=1}^{n-1}\binom{n-1}{\nu}\a^{\nu}\,\Bbb E\bigl[\sigma^{-\nu}_{2k-1}(\bold X)\bigr]
=O\!\left(\sum_{\nu=1}^{n-1}\binom{n-1}{\nu}\left(\a\sqrt{\frac{c}{k}}\right)^{n-1}\right)\\
=O\left[\left(1+\a\sqrt{\frac{c}{k}}\right)^{n-1} - 1\right]=O(nk^{-1/2}).
\end{multline*}
\end{proof}
%It follows that the integral of the big Oh term in \eqref{n1.11} over $C_{2k-1}$ is at most of order
%\begin{equation}\label{n1.31}
%\sum_{\nu=1}^{n-1}\binom{n-1}{\nu}\left(\frac{\a(2/\be)^{1/2}}{k^{1/2}}\right)^{\nu}\!\!
%=\!\left(1+\frac{\a(2/\be)^{1/2}}{k^{1/2}}\right)^{n-1}\!\!\! -1 =O\left(\frac{n}{k^{1/2}}\right),
%\end{equation}
%since $n=O(k^2)$. Therefore, combining \eqref{n1.1}, \eqref{n1.11}, and \eqref{n1.31} we have 
%\begin{equation}\label{n1.8}
%\begin{aligned}
%P_{n,k}&=
%n\int\limits_{\bold x\in C_{2k-1}}\!\!\!\!G^{n-1}(y(\bold x))\, d\bold x +O(n^2 k^{-1/2})\\
%&=n\Bbb E[G^{n-1}(y(\bold X))]+O(n^2 k^{-1/2}),
%\end{aligned}
%\end{equation}
%the big-Oh term being $o(1)$, since $k\gg n^4$. Here $y(\bold X)=\frac{k-\sum_a(1-X_a)}{\left(\sum_a X_a(1-X_a)\right)^{1/2}}$, and $X_a$ are the independent $[0,1]$ Uniforms.
Next, we turn to the (hopefully) leading term in \eqref{n1.11}.
\begin{lemma}\label{nlem2} For $n\ll k^{-1/2}$,
\begin{equation*}
\Bbb E\bigl[G^{n-1}(y(\bold X))\bigr]=\!\int_{-\infty}^{\infty}G^{n-1}(y)\frac{e^{-y^2}}{\sqrt{\pi}}\,dy +O(nk^{-1/2}).
%\frac{\sqrt{2}}{n}\int_0^{\infty} ye^{-y^2/2}\bigl[G^{n}(-y)-G^{n}(y)\bigr]\,dy +
%O(n k^{-1/2}).
\end{equation*}
\end{lemma}
\begin{proof}
By the definition, $y(\bold X)=\frac{k-\sum_a(1-X_a)}{\sqrt{\sum_a X_a(1-X_a)}}$.
Now the mean and the variance of the numerator are $1/2$ and $(2k-1)/12$, and $\Bbb E\left[\sum_a X_a(1-X_a)
\right]=(2k-1)/6$. So, intuitively $y(\bold X)$ is asymptotically $\mathcal N(0,1/2)$, the standard normal with zero mean and variance $1/2$.
To prove this crucial fact, introduce 
\begin{equation}\label{n1.89}
y^*(\bold X):=\frac{k-\sum_a(1-X_a)}{\sqrt{\frac{2k-1}{6}}}.
\end{equation}
Then, 
\begin{multline}\label{n1.85}
y(\bold X)-y^*(\bold x)=\left(k-\sum_a(1-X_a)\right)Z_k(\bold X),\\
Z_k(\bold X)=\frac{\sqrt{\frac{2k-1}{6}}-\sqrt{\sum_a X_a(1-X_a)}}{\sqrt{\sum_a X_a(1-X_a)}\sqrt{\frac{2k-1}{6}}}\\
=\frac{\frac{2k-1}{6}-\sum_a X_a(1-X_a)}
{\sqrt{\sum_a X_a(1-X_a)}\cdot \sqrt{\frac{2k-1}{6}}\left(\sqrt{\sum_a X_a(1-X_a)}+\sqrt{\frac{2k-1}{6}}\right)},\\
\end{multline}
so that 
\[
|Z_k(\bold X)|\le \frac{\left|\frac{2k-1}{6}-\sum_a X_a(1-X_a)\right|}
{\sqrt{\sum_a X_a(1-X_a)}\frac{2k-1}{6}}.
\]
Then, by Cauchy-Schwartz inequality,
\begin{multline}\label{n1.9}
\Bbb E[Z_k^2(\bold X)]\le \Bbb E^{1/2}\left[\left(\frac{2k-1}{6}-\sum_aX_a(1-X_a)\right)^4\right]\\
\times E^{1/2}\left[\left(\sum_a X_a(1-X_a)\right)^{-2}\left(\frac{2k-1}{6}\right)^{-4}\right]\\
\end{multline}
Now, $\sum_a X_a(1-X_a)$ is the sum of the {\it bounded\/} independent random variables, with mean $1/6$, and variance equal $(2k-1)\text{Var}(X(1-X)$. Therefore 
\[
\frac{\frac{2k-1}{6}-\sum_aX_a(1-X_a)}{(2k-1)^{1/2}\text{Var}^{1/2}(X(1-X))}
\]
converges, {\it with all the central moments\/}, to $\mathcal N(0,1)$, meaning that the top expectation in \eqref{n1.9}  is of order $k^2$.

This stronger version of the central limit theorem can be proved by using the Laplace transform of the distribution of $X_a(1-X_a)$, rather than the characteristic function, i. e. the Fourier transform of that distribution.

As for the bottom expectation in \eqref{n1.9}, by the inequality \eqref{n1.31}, it is of order $k^{-6}$.
Hence $\Bbb E[Z_k^2(\bold X)]=O(k^{-2})$. So, applying Cauchy-Schwartz inequality to the top RHS in \eqref{n1.85},
we get a crucial estimate
\begin{equation}\label{n1.86}
\begin{aligned}
\Bbb E\bigl[|y(\bold X)-y^*(\bold X)|\bigr]&\le \Bbb E^{1/2}\left[\left(k-\sum_a(1-X_a)\right)^2\right]\cdot 
\Bbb E^{1/2}[Z_k^2(\bold X)]\\
&=O(k^{1/2} k^{-1})=O(k^{-1/2}).
\end{aligned}
\end{equation}
{\bf Note.\/} In the preceding inequalities we could have used a more general H\H older  inequality $\Bbb E[|UV|]
\le \Bbb E^{1/p}[|U|^p]\Bbb E^{1/q}[|U|^q]$, $(1/p+1/q=1)$, picking up the best $p$ and $q$ in each case. We tried, but the $O(k^{-1/2})$ bound did not budge. 

%Armed with \eqref{n1.86}, we return to \eqref{n1.8}. Since $G(z)\le 1$, we obtain
%\begin{multline*}
%n\bigl|\Bbb E[G^{n-1}(y(\bold X))]-\Bbb E[G^{n-1}(y^*(\bold X)\bigr|\\
%\le n(n-1)\Bbb E\bigl[|y(\bold X)-y^*(\bold X)|\bigr]
%=O(n^2 k^{-1/2}).
%\end{multline*}
%So, \eqref{n1.8} becomes
%\begin{equation}\label{n1.87}
%P_{n,k}=n\,\Bbb E[G^{n-1}(y^*(\bold X))]+ O(n^2 k^{-1/2}).
%\end{equation}
To proceed, define %To fully prepare for the next step we will replace $y^*(\bold X)$ defined in \eqref{n1.89} with 
\[
y^{**}(\bold X)=\frac{(2k-1)/2-\sum_a(1-X_a)}{\sqrt{\frac{2k-1}{6}}};
\]
clearly $|y^{**}(\bold X)-y^*(\bold X)|=O(k^{-1/2})$.
%, so that \eqref{n1.86} implies 
%\begin{equation}\label{n1.87}
%\Bbb E\bigl[|y(\bold X)-y^{**}(\bold X)|\bigr]=O(k^{-1/2}). 
%\end{equation}
%as well. 
Since the mean and the variance of $\sum_a(1-X_a)$ are $(2k-1)/2$ and $\frac{2k-1}{12}$, by the central limit theorem $y^{**}(\bold X)$ is asymptotically $\mathcal N(0,1/2)$. In other words, 
%since 
%\begin{equation}\label{n1.875}
%n\Bbb E[G^{n-1}(y^*(\bold X))]-n\Bbb E[G^{n-1}(y^{**}(\bold X))]=O(n^2k^{-1/2}).
%\end{equation}
introducing $F_{y{**}(\bold X)}(y)$, the distribution function of $y^{**}(\bold X)$, we have
\[
\lim_{k\to\infty}F_{y{**}(\bold X)}(y)=\Bbb P(\mathcal N(0,1/2)\le y)=\frac{1}{\sqrt{\pi}}\int_{-\infty}^ye^{-z^2}\,dz=:\Psi(y).
\]
Moreover, since $\Bbb E[|X-1/2|^3]<\infty$, by a general theorem (Feller \cite{Fel}, Ch. XVI, Section 4): uniformly for all $y$,
\[
F_{y^{**}(\bold X)}(y)=\Psi(y) +\frac{c\mu_3}{\sigma^3k^{1/2}}(1-2y^2)e^{-y^2}+o(k^{-1/2}),
\]
where $c$ is an absolute constant, and $\sigma^2$ and $\mu_3$ are the second and the third central moments of $1-X_a$, meaning that in our case, $\sigma^2=1/6$, and $\mu_3=0$, whence $F_{y{**}(\bold X)}(y)=\Psi(y)+o(k^{-1/2})$. In fact, since the characteristic function $\Bbb E\bigl[e^{iu(X-1/2)}\bigr]$ is square integrable, the density $p_{y^{**}}(y)$ of $y^{**}(\bold X)$ is given by a similar asymptotic formula,
(\cite{Fel}, Ch. XVI, Section 2), which for the case $\mu_3=0$ simplifies to 
$p_{y^{**}}(y)%=\frac{e^{-y^2/2}}{\sqrt{2\pi}}-\frac{\mu_3}{6\sigma^3 k^{1/2}}(y^3-3y)\frac{e^{-y^2/2}}{\sqrt{2\pi}}+o(k^{-1/2})\\
=\pi^{-1/2} e^{-y^2}+o(k^{-1/2})$.
In particular, the density $p_{y^{**}}(y)$ is uniformly bounded for all $k$ and $y$. 

So, recalling that $|y^{**}(\bold X)-y^*(\bold X)|\le c k^{-1/2}$, we have 
\begin{multline*}
F_{y^*(\bold X)}(y)=\Bbb P(y^*(\bold X)\le y)\le F_{y^{**}(\bold X)}(y+c k^{-1/2})\\
=F_{y^{**}(\bold X)}(y)+\int_y^{y+c k^{-1/2}}p_{y^{**}(\bold X)}(z)\,dz)\\
=F_{y^{**}(\bold X)}(y)+O(k^{-1/2}).
\end{multline*}
Likewise, $F_{y^*(\bold X)}(y)\ge F_{y^{**}(\bold X)}(y)+O(k^{-1/2})$, so that $F_{y^*(\bold X)}(y)= F_{y^{**}(\bold X)}(y)+O(k^{-1/2})$, enabling us to conclude that 
$F_{y^*(\bold X)}(y)=\Psi(y)+O(k^{-1/2})$, 
%Since  $F_{y^{**}(\bold X)}(y)=F(y)+o(k^{-1/2})$, we conclude that
%$F_{y^*(\bold X)}(y)=F(y)+O(k^{-1/2})$ 
uniformly for all $y$ and $k$.
Therefore 
\begin{multline*}
\Bbb E\bigl[G^{n-1}(y^{*}(\bold X))\bigr]=\int_{-\infty}^{\infty}G^{n-1}(y)\,dF_{y^{*}(\bold X)}(y)\\
=G^{n-1}(y)F_{y^{*}(\bold X)}(y)\big|_{-\infty}^{\infty}-\int_{-\infty}^{\infty}(G^{n-1}(y))'F_{y^{*}(\bold X)}(y)\,dy\\
=-\int_{-\infty}^{\infty}\bigl[(G^{n-1}(y))' (\Psi(y)+O(k^{-1/2}))\bigr]\,dy\\
=-\int_{-\infty}^{\infty}(G^{n-1}(y))' \Psi(y)\,dy +O\bigl[k^{-1/2}\bigl(G^{n-1}(-\infty)-G^{n-1}(\infty)\bigr)\bigr]\\
=\int_{-\infty}^{\infty}G^{n-1}(y)\frac{e^{-y^2}}{\sqrt{\pi}}\,dy+O(k^{-1/2}).
%\\
%=\frac{G^n(y)}{n}\big|_{\infty}^{-\infty}+O(k^{-1/2})=\frac{1}{n}+O(k^{-1/2}).
\end{multline*} 
Finally, 
\begin{multline*}
\big |G^{n-1}(y(\bold X))-G^{n-1}(y^*(\bold X))\big|\le (n-1)\big|G(y(\bold X))-G(y^*(\bold X))\big|\\
=O\bigl(n|y(\bold X)-y^*(\bold X)|).
\end{multline*}
So, using $\Bbb E\bigl[|y^*(\bold X)-y(\bold X)|\bigr]=O(k^{-1/2})$, we obtain 
\begin{equation*}
\Bbb E\bigl[G^{n-1}(y(\bold X))\bigr]\!=\Bbb E\bigl[G^{n-1}(y^*(\bold X))\bigr]+O(nk^{-1/2})\!=\!I_{n}+O(\nu k^{-1/2}),
\end{equation*} 
$I_{n}:=\!\int_{-\infty}^{\infty}G^{n-1}(y)\frac{e^{-y^2}}{\sqrt{\pi}}\,dy$ . 
\end{proof}
It remains to bound $I_n$.
\begin{lemma}\label{nlem3} $I_n=O\bigl(\exp(-n^{\eps(n)})\bigr)$ for every $\eps(n)\downarrow 0$. In particular, $I_n$ is super-polynomially small.
\end{lemma}
\begin{proof} 
First,  
\begin{multline*}
%\int_{-\infty}^{\infty}G^{\nu}(y)\frac{e^{-y^2}}{\sqrt{\pi}}\,dy
I_{n}=\frac{\sqrt{2}}{n}
%\int_{-\infty}^{\infty}n G^{n-1}(y)\frac{e^{-y^2}/2}{\sqrt{2\pi}}\cdot e^{-y^2/2}\,dy\\
\int_{-\infty}^{\infty}(-G^n(y))' e^{-y^2/2}\,dy\\
=\frac{\sqrt{2}}{n}\int_{-\infty}^{\infty}(-y)e^{-y^2/2}G^{n}(y)\,dy\\
=\frac{\sqrt{2}}{n}\int_0^{\infty} ye^{-y^2/2}\bigl[G^{n}(-y)-G^{n}(y)\bigr]\,dy.
\end{multline*}
%First of all, $I_{\nu}=o(\nu^{-1})$, since $\int_0^{\infty}ye^{-y^2/2}\,dy=1$, and $G^{\nu+1}(-y)\to 0$ for every fixed $y$. More precisely,
Since $G(y)\le 1/2$ for $y\ge 0$, we have $
\int_0^{\infty} ye^{-y^2/2}G^{n}(y)\,dy=O(2^{-n})$. So $I_n=\frac{\sqrt{2}}{n}(J_n+O(2^{-n}))$, where 
\[
J_n:=\int_0^{\infty}ye^{-y^2/2} G^{n}(-y)\,dy=\int_0^{\infty}ye^{-y^2/2}\left(1-\frac{1}{\sqrt{2\pi}}\int_y^{\infty}e^{-z^2/2}\,dz\right)^n\,dy.
\]
Because of how the integrand's depends on $n$, we expect that the dominant contribution to $J_n$ comes from a relatively small neighborhood of the integrand's maximum point. It is easy to see that the integrand is strictly log-concave, and since the integrand approaches zero both at $y=0$ and $y=\infty$, it attains its unique maximum at a unique point. A closer look shows that the maximum
is very near $\bar y:=\left(2\log\left(\frac{n}{\sqrt{\log n}}\right)\right)^{1/2}$. We do not have to prove this clam, since its sole purpose is to 
{\it motivate\/} a rigorous argument. 
Let $J_{n,1}$ and $J_{n,2}$ denote the contributions to $J_n$ coming respectively from $y\in [0,y(\a)]$ and $y\in [y(\a),\infty)$, $y(\a):=\a \bar y$. Here $\a=\a(n)\uparrow 1$ will be chosen based on bounds on $J_{n,j}$.

{\bf (a)\/} Using a classic inequality
 \begin{equation}\label{verynew1}
\int_y^{\infty}e^{-z^2/2}\,dz\ge \bigl(y^{-1}-y^{-3}\big)e^{-y^2/2},
\end{equation}
we have: 
\begin{multline}\label{new4}
J_{n,1}\le\left(1-\frac{1}{\sqrt{2\pi}}\int_{y(\a)}^{\infty}e^{-z^2/2}\,dz\right)^n\int_0^{y(\a)}ye^{-y^2/2}\,dy\\
\le \exp\left(-\frac{n}{\sqrt{2\pi}}\int_{y(\a)}^{\infty}e^{-z^2/2}\,dz\right)\\
\le\exp\left(-\bigl(y^{-1}(\a)-y^{-3}(\a)\bigr)\frac{n}{\sqrt{2\pi}}e^{-y^2(\a)/2}\right).
\end{multline}
%\int_2^{y(\eps)} y\exp\left(-\frac{y^2}{2}-\frac{n}{\sqrt{2\pi}}\left(\frac{1}{y}-\frac{1}{y^3}\right)e^{-y^2/2}\right)\,dy\\
%\le \int_0^{y(\eps)}\!\! y\exp\left(\!-\frac{y^2}{2}-\frac{n}{\sqrt{2\pi}}\frac{3}{4y(\eps)}e^{-y^2/2}\!\right)\,dy\\
%=\int_0^{y^2(\eps)/2}\!\!\!\exp\left(\!-z-\frac{n}{\sqrt{2\pi}}\frac{3}{4y(\eps)}e^{-z}\!\right) dz\\
%=\int_{e^{-y^2(\eps)/2}}^1\exp\left(-\frac{n}{\sqrt{2\pi}}\frac{3}{4y(\eps)}\eta\right)\,d\eta\\
%\le \frac{\exp\left(-\frac{n}{\sqrt{2\pi}}\frac{3}{4y(\eps)}e^{-y^2(\eps)/2}\right)}{\frac{n}{\sqrt{2\pi}}\frac{3}{4y(\eps)}}.
%\end{multline*}
%And obviously the contribution to $J_{n,1}$ from $[0,2]$ is at most $q^n$, where $q:=1-\frac{1}{\sqrt{2\pi}}\int_2^{\infty}e^{-z^2/2}\,dz$.

{\bf (b)\/} Next,
\begin{multline*}
J_{n,2}=\int_{y(\a)}^{\infty}ye^{-y^2/2}\left(1-\frac{1}{\sqrt{2\pi}}\int_y^{\infty}e^{-z^2/2}\,dz\right)^n\,dy\\
\le \int_{y(\a)}^{\infty}ye^{-y^2/2}
\exp\left(-\frac{n}{\sqrt{2\pi}}\int_y^{\infty}e^{-z^2/2}\,dz\right)\,dy\\
=\int_{y^2(\a)/2}^{\infty}\exp\left(-u-\frac{n}{\sqrt{2\pi}}\int_{\sqrt{2u}}^{\infty}e^{-z^2/2}\,dz\right)\,du\\
=\int_{0}^{e^{-y^2(\a)/2}}e^{-nK(v)}\,dv, \quad K(v):=(2\pi)^{-1/2}\int_{\sqrt{2\log(1/v)}}^{\infty}e^{-z^2/2}\,dz.
\end{multline*}
Now, $K'(v)=\frac{1}{2\sqrt{\pi}}\log^{-1/2}(1/v)(>0)$ increases with $v$, so that $K(v)$ is convex; hence
$
K(v)\ge K(e^{-y^2(\a)/2})+K'(e^{-y^2(a)/2})\bigl(v-e^{-y^2(\a)/2}\bigr).
$
Therefore 
$
J_{n,2}\le \frac{e^{-nK(e^{-y^2(\a)/2})}}{nK'(e^{-y^2(\a)/2})}.
$
Here
\begin{align*}
nK(e^{-y^2(\a)/2})&=\frac{n}{\sqrt{2\pi}}\int_{y(\a)}^{\infty}e^{-z^2/2}\,dz\\
&\ge \bigl(y^{-1}(\a) -y^{-3}(\a)\bigr)\frac{n}{\sqrt{2\pi}}e^{-y^2(\a)/2},\\
%\ge c_2\frac{n}{\sqrt{\log n}}\cdot\exp\Bigl(-(1-\eps)^2\log\bigl(\frac{n}{\sqrt{\log n}}\bigr)\Bigr)=c_2\Bigl(\frac{n}{\sqrt{\log n}}\Bigr)^{2\eps -\eps^2}\to\infty,\\
nK'(e^{-y^2(\eps)/2})&=\frac{n}{\sqrt{2\pi}}y^{-1}(a).
\end{align*}
So, $J_{n,2}$ is at most of order
\begin{equation}\label{new5}
\frac{\exp\left(-\bigl(y^{-1}(\a) -y^{-3}(\a)\bigr)\frac{n}{\sqrt{2\pi}}e^{-y^2(\a)/2}\right)}{ny^{-1}(\a)}
\end{equation}
Here $ny^{-1}(\a)$ is exactly of order $n\log^{-1/2} n\to \infty$. So \eqref{new4} and \eqref{new5} together imply that 
\begin{multline*}
J_n=O\left[\exp\left(-\bigl(y^{-1}(\a)-y^{-3}(\a)\bigr)\frac{n}{\sqrt{2\pi}}e^{-y^2(\a)/2}\right)\right]\\
=O\left[\exp\left(-\frac{1+o(1))}{2\sqrt{\pi}}\frac{n}{\sqrt{\log n}}\left(\frac{\sqrt{\log n}}{n}\right)^{\a^2}\right)\right]\\
=O\left[\exp\left(-\frac{1+o(1))}{2\sqrt{\pi}}\left(\frac{n}{\sqrt{\log n}}\right)^{1-\a^2}\right)\right]
\end{multline*}
%\le \exp\lbigl(-\exp(\log^{-s}n)\bigr),
For $\a:=1- \eps(n)$, $\eps(n)\downarrow 0$, we have $1-\a^2=2\eps(n)-\eps^2(n)$, whence $J_n\le \exp\bigl(-n^{\eps(n)}\bigr)$.
We conclude that
\[
I_n=\frac{\sqrt{2}}{n}\bigl(J_n+O(2^{-n})\bigr)=O\bigl[\exp\bigl(-n^{\eps(n)})\bigr].
\]
In particular, selecting $\eps(n)=\frac{\log(\log(\ell\log n))}{\log n}$, $\ell>0$ being fixed. we obtain $I_n=O(n^{-\ell})$. This completes the proof of Lemma \ref{nlem2}.
\end{proof}
Combining Lemma \ref{nlem1} and Lemma \ref{nlem2}, and assuming $k\gg n^2$, we transform \eqref{n1.11} into
\begin{multline*}
n^{-1}Q_{2k-1,n}=\int_{-\infty}^{\infty}\!\!G^{n-1}(y)\frac{e^{-y^2}}{\sqrt{\pi}}\,dy+O(nk^{-1/2})\\
=O\bigl[\exp\bigl(-n^{\eps(n)}\bigr)+nk^{-1/2}\bigr].
\end{multline*}
%+O(nk^{-1/2})+O\!\left(\sum_{\nu=1}^{n-1}\binom{n-1}{\nu}\left(\a\sqrt{\frac{c}{k}}\right)^{\nu}\right)\\
%=\int_{-\infty}^{\infty}G^{n-1}(y)\frac{e^{-y^2}}{\sqrt{\pi}}\,dy+O(nk^{-1/2}).
%\end{multline*}
Thus, for $m=2k-1$, and $m\gg n^4$,
\begin{multline*}
Q_{m,n}\!=\!n\!\!\int\limits_{-\infty}^{\infty}\!\!\!G^{n-1}(y)\frac{e^{-y^2}}{\sqrt{\pi}}\,dy+O(n^2m^{-1/2})\\
=O\bigl[\exp\bigl(-n^{\eps(n)}\bigr)+n^2 m^{-1/2}\bigr].
%\quad G(y):=\frac{1}{\sqrt{2\pi}}\int_y^{\infty}\!\!e^{-z^2/2}\,dz.
\end{multline*}
And this formula holds for even $m$ as well. This completes the proof of the theorem.

\noindent {\bf Acknowledgment.\/} The author learned about the Condorcet winner problem from \cite{Sau} posted
by Lisa Sauermann on math arXiv in 2022. I owe debt of gratitude to \'Elie de Panafieu for catching a consequential error in an original manuscript, and for a very valuable feedback overall. I thank Ilia Tsetlin for a cautionary comment on the original claim.

\end{document}